\documentclass[10pt]{amsart}
\usepackage{latexsym, amsmath,amssymb}

\renewcommand{\epsilon}{\varepsilon}
\newcommand{\newsection}[1]
{\subsection{#1}\setcounter{theorem}{0} \setcounter{equation}{0}
\par\noindent}

\newtheorem{theorem}{Theorem}

\newtheorem{lemma}[theorem]{Lemma}
\newtheorem{corr}[theorem]{Corollary}

\newtheorem{proposition}[theorem]{Proposition}
\newtheorem{deff}[theorem]{Definition}

\newcommand{\bth}{\begin{theorem}}
\newcommand{\ble}{\begin{lemma}}
\newcommand{\bcor}{\begin{corr}}

\newcommand{\bdeff}{\begin{deff}}

\newcommand{\bprop}{\begin{proposition}}
\newcommand{\ele}{\end{lemma}}
\newcommand{\ecor}{\end{corr}}
\newcommand{\edeff}{\end{deff}}

\newcommand{\eprop}{\end{proposition}}

\newcommand{\la}{\lambda}

\newcommand{\supp}{\text{supp }}
\renewcommand{\Pi}{\varPi}

\renewcommand{\epsilon}{\varepsilon}

\newcommand{\R}{{\mathbb R}}
\newcommand{\Z}{{\mathbb Z}}
\newcommand{\N}{{\mathbb N}}

\begin{document}

\title[Localized $L^p$-estimates of eigenfunctions]
{Localized $L^p$-estimates of eigenfunctions: \\ A note on an article of Hezari and Rivi\`ere}
\thanks{The author was supported in part by the NSF grant DMS-1361476.}
%
%
%
%
%
%
\author{Christopher D. Sogge}
\address{Department of Mathematics,  Johns Hopkins University,
Baltimore, MD 21218}
\email{sogge@jhu.edu}
\keywords{Eigenfunctions, Quantum ergodicity, Negative curvature}
\subjclass[2010]{Primary 58J51; Secondary 35A99, 42B37}

\maketitle

\begin{abstract}
	We use a straightforward variation on a recent argument of Hezari and Rivi\`ere~\cite{HR} to obtain localized $L^p$-estimates for all exponents larger
	than or equal to the  critical exponent 
	$p_c=\tfrac{2(n+1)}{n-1}$.  We are able to this  directly by just  using   the $L^{p}$-bounds for spectral projection operators from our much earlier work \cite{Seig}.
 	The localized bounds we obtain here imply, for instance, that, for a density one sequence of eigenvalues on a
	manifold whose geodesic flow is ergodic,  all of the $L^p$, $2<p\le \infty$, bounds 
	of the corresponding eigenfunctions are relatively small compared to the general ones in
	\cite{Seig}, which are saturated on round spheres.  The connection with quantum ergodicity was established for exponents  $2<p<p_c$ in the recent
	results of the author \cite{SK} and Blair and the author \cite{BS2}; however, the article of Hezari and Rivi\`ere~\cite{HR} was the first one to make this connection (in the case of negatively curved manifolds) for the critical exponent, $p_c$. As is well known, and we indicate here, bounds for the critical exponent, $p_c$, imply ones for all of the other exponents $2<p\le \infty$.  The localized estimates involve $L^2$-norms over
	small geodesic balls $B_r$ of radius $r$, and we shall go over what
	happens for these in certain model cases on the sphere and on 
	manifolds of nonpositive curvature.  We shall also state a problem
	as to when one can improve on the trivial $O(r^{\frac12})$ estimates for these
	$L^2(B_r)$ bounds.  If $r=\la^{-1}$, one can improve on the trivial estimates if one has improved
	$L^{p_c}(M)$ bounds just by using H\"older's inequality; however, obtaining improved bounds for $r\gg \la^{-1}$ seems to be subtle.
\end{abstract}

\newsection{Introduction}

Let $(M,g)$ be an $n$-dimensional compact manifold without boundary with
$n\ge 2$.  Then if $\Delta_g$ is the associated Laplace-Beltrami
operator, we shall consider $L^2$-normalized eigenfunctions of
$\sqrt{-\Delta_g}$, i.e., functions $e_\la$ satisfying
\begin{equation}\label{1.1}
-\Delta_g e_\la(x)=\la^2 e_\la(x), \quad\text{and } \, \,
\int_M |e_\la|^2 \, dV_g=1.
\end{equation}
Here $dV_g$ denotes the volume element for $(M,g)$, and, in what
follows, all of the $L^p$-norms are taken with respect to this measure.

Our main result says that one can control the critical
$L^p$-norms of eigenfunctions in terms of local $L^2$-estimates over
balls of possibly small size.

\begin{theorem}\label{loctheorem}
  For $r>0$ smaller than the injectivity radius of $(M,g)$,
let $B_r(x)$ denote the geodesic ball of radius $r$ centered at $x$.
Then	there is a uniform constant $C$, depending only on $(M,g)$, so
that for $\la \ge 1$ and eigenfunctions as in \eqref{1.1} we have
\begin{equation} \label{1.2}
	\|e_\la\|_{L^{\frac{2(n+1)}{n-1}}(M)}\le C \la^{\frac{n-1}{2(n+1)}} \, 
	\Bigl(r^{-\frac{n+1}4}\sup_{x\in M} \|e_\la\|_{L^2(B_r(x))}
	\Bigr)^{\frac2{n+1}}, \quad
	\la^{-1}\le r\le \text{Inj }M,
\end{equation}
where $\text{Inj }M$ denotes the injectivity radius of $(M,g)$.
\end{theorem}

The special case of \eqref{1.2} corresponding to $r\approx 1$ is
equivalent to the earlier estimates of the author \cite{Seig},
\begin{equation}\label{1.3}
	\|e_\la\|_{L^{\frac{2(n+1)}{n-1}}(M)}\le C\la^{\frac{n-1}{2(n+1)}}, 
	\qquad \la \ge 1,
\end{equation}
which are saturated on round spheres both by zonal spherical harmonics and highest weight spherical harmonics.  Note that by a Bernstein inequality,
\eqref{1.3} yields the sup-norm estimates
$$\|e_\la\|_{L^\infty(M)}\le C\la^{\frac{n-1}{2(n+1)}}\la^{\frac{n(n-1)}{2(n+1)}}=C\la^{\frac{n-1}2},$$
and so by interpolating between this estimate, the trivial $L^2$
estimate and \eqref{1.3}, we obtain the results of \cite{Seig}:
\begin{equation}
	\label{1.4}
	\|e_\la\|_{L^p(M)}\le C\la^{\sigma(p)},
	\quad \la \ge 1, \end{equation}
where
\begin{equation}
	\label{1.5}
	\sigma(p)=\begin{cases}
	n(\tfrac12-\tfrac1p)-\tfrac12, \quad  \tfrac{2(n+1)}{n-1}\le p\le \infty,
	\\ \\
	\tfrac{n-1}2(\tfrac12-\tfrac1p), \qquad  2\le p\le \tfrac{2(n+1)}{n-1}.	
	\end{cases}
\end{equation}
As was shown in \cite{Sthesis}, these estimates are also saturated on the
round sphere.  To be more specific, for $2<p\le \tfrac{2(n+1)}{n-1}$ they are saturated by
the highest weight spherical harmonics, while for 
$\tfrac{2(n+1)}{n-1}\le p\le \infty$, they are saturated by zonal spherical
harmonics.

In \cite{Seig}, a stronger version of \eqref{1.4}--\eqref{1.5} was obtained for general
$(M,g)$.
Specifically, if $E_j$ denotes the projection onto the eigenspace
of $\sqrt{-\Delta_g}$
with eigenvalue $\la_j$, and if
$\chi_\la$ denotes the spectral projection operator,
$$\chi_\la f=\sum_{\la_j\in [\la,\la+1)}E_jf, \quad \la \ge 0,$$
projecting onto unit bands of frequencies, then it was shown in \cite{Seig}
that 
\begin{equation}\label{1.6}
	\|\chi_\la f\|_{L^p(M)}\le C(1+\la)^{\sigma(p)}\|f\|_{L^2(M)},
\end{equation}
if $\sigma(p)$ is as in \eqref{1.5}.  Here
$$0=\la_0<\la_1\le \la_2\le \cdots,$$
denotes the spectrum of $\sqrt{-\Delta_g}$ labeled with respect to multiplicity, to which we can associate an orthonormal basis of eigenfunctions
$\{e_{\la_j}\}_{j=0}^\infty$.

By the argument that we just gave
showing how \eqref{1.4} follows from \eqref{1.3}, the preceding estimates
just follow from the special case
\begin{equation}
	\label{1.7}
	\|\chi_\la f\|_{L^{\frac{2(n+1)}{n-1}}(M)}\le C
	(1+\la)^{\frac{n-1}{2(n+1)}}\|f\|_{L^2(M)}.
\end{equation}
We shall prove the localized estimates \eqref{1.2} for eigenfunctions
just by using \eqref{1.7}, and, thus, unlike the arguments in \cite{HR},
avoid the use of semi-classical analysis.

Before doing this, let us record a corollary of \eqref{1.2}.  

\begin{corr}\label{corr}
Assume that the geodesic flow on $(M,g)$  is ergodic.
Then there exists a density one subsequence of eigenvalues $\la_{j_k}$
so that for every $2<p\le \infty$ we have
\begin{equation}
	\label{1.8}
	\|e_{\la_{j_k}}\|_{L^p(M)}=o(\la_{j_k}^{\sigma(p)}),
\end{equation}
if $\sigma(p)$ is as in \eqref{1.5}.	
\end{corr}

To see this, we note that, by the argument that we just gave, \eqref{1.8} follows from the special case
\begin{equation}
	\label{1.9}
	\|e_{\la_{j_k}}\|_{L^{\frac{2(n+1)}{n-1}}(M)}=o(\la_{j_k}^{\frac{n-1}{2(n+1)}}).
\end{equation}
To prove this, we use the quantum ergodicity theorem of Colin de Verdi\`ere--Shnirelman--Zelditch~\cite{Colin}--\cite{Snirelman}--\cite{Zelditch} to select a density one subsequence of eigenvalues
so that the corresponding eigenfunctions satisfy (see e.g., 
\cite[Corollary 6.2.4]{SoHang})
\begin{equation}\label{1.10}\int_\Omega |e_{\la_{j_k}}|^2 \, dV_g \to |\Omega|/|M|, \quad k\to \infty,
\end{equation}
for Jordan measurable subsets $\Omega$ of $M$, where $|\Omega|$ denotes
its $dV_g$-measure.  Since $|B_r(x)|\approx r^n$ with bounds independent of
$x\in M$ if $r\ll 1$, we get that
$$\lim_{k\to \infty}r^{-{\frac{n+1}4}}\|e_{\la_{j_k}}\|_{L^2(B_r(x))}
\approx r^{\frac{n-1}4}.$$
This along with \eqref{1.2} yields \eqref{1.9} since for
a given fixed $r\ll 1$ we can find $O(r^{-n})$ points $x_\ell\in M$
so  that the resulting balls $B_r(x_\ell)$ cover $M$ and have
overlap of at most a constant $N=N((M,g))$, which can be chosen independent
of $r$.

As we shall see in \S 4, for all $(M,g)$ there is the trivial uniform bounds
\begin{equation}\label{1.11}
\|e_\la\|_{L^2(B_r(x))}\le Cr^{\frac12},
\quad \la^{-1}\le r\le \text{Inj }M, \, \, x\in M.
\end{equation}
As we shall also show, this bound is saturated by the $L^2$-normalized
zonal functions on $S^n$, $Z_\la$, $\la=\sqrt{k(k+n-1)}$, $k\in \N$,
centered at a given $x_0\in S^n$, since
$$
\|Z_\la\|_{L^2(B_r(x_0))}\approx r^{\frac12}, \quad
\la^{-1}\le r\le \text{Inj }M.
$$
On the other hand, in the case of negative curvature, recent
results of Han~\cite{H} and Hezari and Rivi\`ere~\cite{HR} improve
upon \eqref{1.11} considerably in the sense that, if all the 
sectional curvatures of $(M,g)$ are negative then there is always
a density one sequence of eigenvalues $\{\la_{j_k}\}$ such that
one has the small-scale quantum ergodic estimates
$$
\|e_{\la_{j_k}}\|_{L^2(B_r(x))}\approx r^n, \, \, \,
\forall x\in M, \, \, \text{if } \, \, r=(\log \la_{j_k})^{-\kappa},
$$
for a range of powers $\kappa>0$ depending on the dimension.  
Also, B\'erard's~\cite{PB} proof of improved error estimates for
the Weyl formula for manifolds of nonpositive curvature imply that
one can always improve \eqref{1.11} for the smallest allowable
$r$ if $(M,g)$ has nonpositive curvature, since then one has
$$
\|e_\la\|_{L^2(B_r(x))}\le Cr^{\frac12}/(\log \la )^{\frac12},
\quad \text{if } \, \, r=\la^{-1}.
$$

It would be interesting to find other general cases where \eqref{1.11}
can be improved.

We note just by using H\"older's inequality, if one has the improved estimates \eqref{1.9} for the critical exponent $p_c=\tfrac{2(n+1)}{n-1}$, then one automatically has an improvement over \eqref{1.11} at the smallest possible scale, i.e., 
$$\|e_{\la_{j_k}}\|_{L^2(B_{\la_{j_k}^{-1}}(x))}=o(\la_{j_k}^{-\frac12}), \quad
x\in M.$$
The converse assertion need not hold, though.  For the $L^2$-normalized highest weight spherical harmonics satisfy 
$\|Q_\la\|_{L^2(S^n)}\approx \la^{\frac{n-1}{2(n+1)}}$ as well as $\|Q_\la\|_{L^2(B_r(x))}=o(r^{\frac12})$ for all $\la^{-1}\le r\ll \la^{-\frac12}$ (see
\S 4).

\newsection{Proof of the localized $L^{p_c}$-bounds}

Choose a nonnegative function $\rho\in {\mathcal S}(\R)$ satisfying
\begin{equation}
	\label{2.1}
	\rho(0)=1 \quad \text{and } \, \, 
	\Hat \rho(t)=0, \, \, \text{if } \, |t|\ge 1.
\end{equation}
It then follows that if we let $P=\sqrt{-\Delta_g}$ then the operator
defined by
\begin{equation}\label{2.2}T_{\la,r} = \frac1\pi \int_{-\infty}^\infty r^{-1} \Hat\rho(r^{-1}t) \, e^{it\la}
\cos(tP) \, dt,\end{equation}
by Euler's formula equals $\rho(r(\la-P))+\rho((r(\la+P))$.  Therefore,
by \eqref{2.1}
$$T_{\la,r} e_\la =\bigl[1+\rho(2r\la)\bigr] \, e_\la,$$
and since we are assuming that $\rho$ is nonnegative we have
$|T_{\la,r} e_\la|\ge |e_\la|$, and so
\begin{equation}
	\label{2.3}
	\|e_\la \|_{L^{\frac{2(n+1)}{n-1}}(B_r(x))}\le \|T_{\la,r} e_\la\|_{L^{\frac{2(n+1)}{n-1}}(B_r(x))}.
\end{equation}

Next, we note that by Huygens' principle, the kernel $\bigl(\cos tP\bigr)(x,y)$
vanishes if the geodesic distance between $x$ and $y$ is greater than $t$.
Therefore, we conclude from the second part of \eqref{2.1} that the kernel
$K_{\la,r}(x,y)$ of $T_{\la,r}$ satisfies
\begin{equation}\label{2.4}K_{\la,r}(x,y)=0, \quad \text{if } \, \, d_g(x,y)>r,	
\end{equation}
where $d_g$ denotes geodesic distance in $(M,g)$.  Consequently, if we could
show that there is a uniform constant $C$ so that when $\la \ge 1$
and $\la^{-1}\le r\le \text{Inj }M$ 
\begin{equation}
	\label{2.5}
	\|T_{\la,r}f\|_{L^{\frac{2(n+1)}{n-1}}(M)}\le Cr^{-\frac12}\la^{\frac{n-1}{2(n+1)}}\|f\|_{L^2(M)},
\end{equation}
then it follows from \eqref{2.3} and \eqref{2.5} that we would have the uniform localized
estimates
\begin{equation}
	\label{2.6}
	\|e_\la \|_{L^{\frac{2(n+1)}{n-1}}(B_r(x))}\le Cr^{-\frac12}\la^{\frac{n-1}{2(n+1)}}\|e_\la\|_{L^2(B_{2r}(x))}, 
	\quad \la \ge 1, \, \, \, \la^{-1}\le r\le \text{Inj }M.
\end{equation}

Let us postpone the proof of \eqref{2.5} for the moment, which will be
a simple consequence of the author's earlier estimate \eqref{1.7}, and see, now,
how \eqref{2.6} implies our main estimate \eqref{1.2}.  Clearly, it suffices
to prove the latter when $\la^{-1}\le r\le \delta$ where $\delta$ is a fixed
positive number since the bounds for $\delta <r<\text{Inj }M$ follow from 
\eqref{1.3}.  To do this we use the fact that if $r$ is small enough we can
cover $M$ by balls $\{B_r(x_\ell)\}_{\ell=1}^{N(r)}$ where
$N(r)\approx r^{-n}$ and where the doubled balls $\{B_{2r}(x_\ell)\}_{\ell=1}^{N(r)}$ overlap at most $A$ times, with $A$ being a constant that depends
on $(M,g)$ but not on small $r>0$.  We then conclude from \eqref{2.6} that,
if $C_0=C^{\frac{2(n+1)}{n-1}}$, then
\begin{align*}
\|e_\la\|_{L^{\frac{2(n+1)}{n-1}}(M)}^{	\frac{2(n+1)}{n-1}}
&\le \sum_{\ell =1}^{N(r)}\|e_\la\|_{L^{\frac{2(n+1)}{n-1}}(B_r(x_\ell))}^{	\frac{2(n+1)}{n-1}}
\\
&\le C_0 \la \,  r^{-\frac{n+1}{n-1}}\sum_{\ell =1}^{N(r)}\|e_\la\|_{L^2(B_{2r}(x_\ell))}^{\frac{2(n+1)}{n-1}}
\\
&\le C_0 \la \,  r^{-\frac{n+1}{n-1}}\,  \Bigl(\sup_{1\le l\le N(r)}
\|e_\la\|_{L^2(B_{2r}(x_\ell))}^{\frac{2(n+1)}{n-1}-2}\Bigr)\sum_{\ell =1}^{N(r)}
\|e_\la\|_{L^2(B_{2r}(x_\ell))}^2
\\
&\le AC_0 \la \,  r^{-\frac{n+1}{n-1}} \, \Bigl(\sup_{1\le l\le N(r)}
\|e_\la\|_{L^2(B_{2r}(x_\ell))}^{\frac4{n-1}}\Bigr)\, 
\|e_\la\|_{L^2(M)}^2,
\end{align*}
which of course yields \eqref{1.2} because of \eqref{1.1}.

\begin{proof}[End of Proof of Theorem~\ref{loctheorem}]

To complete the proof of 
our main result we just need to prove \eqref{2.5}.

We first recall that
\begin{equation}
	\label{2.7}
	T_{\la,r}f=\sum_{j=0}^\infty \bigl[\rho(r(\la-\la_j))+\rho(r(\la+\la_j))\bigr] \, E_jf.
\end{equation}
Since $\rho\in {\mathcal S}(\R)$, for every $N=1,2,3,\dots$, we have for $\la\ge1$
\begin{equation}
	\label{2.8}
	|\rho(r(\la-\la_j))|+|\rho(r(\la+\la_j))|\le C_N(1+r|\la-\la_j|)^{-N}.
\end{equation}
Therefore,
\begin{equation}
	\label{2.9}
	\|\chi_kT_{\la,r}f\|_{L^2(M)}\le C_N(1+r|\la-k|)^{-N}\|\chi_k f\|_{L^2(M)}, \, \, \, N=1,2,\dots .
\end{equation}

To exploit this, let, for $\ell \in \Z$,
$$I_\ell =\Bigl[\la +r^{-1}(2\ell-1),\la+r^{-1}(2\ell + 1)\Bigr).$$
Then $\R =\bigcup_{\ell}I_\ell$.  Also, since $O(r^{-1})$ intervals $[k-1,k)$ intersect a given $I_\ell$ as $k$ ranges over $\N$,
we can use
Minkowski's inequality and the Cauchy-Schwarz inequality to see that 
$$
\Bigl\|\, \sum_{\{k: \, [k-1,k)\cap I_\ell \ne \emptyset\}}\chi_k h
\Bigr\|_{L^{\frac{2(n+1)}{n-1}}(M)}
\le Cr^{-\frac12}\Bigl( \, \sum_{\{k: \, [k-1,k)\cap I_\ell \ne \emptyset\}} \|\chi_k h\|_{L^{\frac{2(n+1)}{n-1}}(M)}^2
\, \Bigr)^{\frac12}.
$$
If we use this along with \eqref{1.6} and \eqref{2.8}, we deduce that
\begin{align*}
	\label{2.10}
	\Bigl\|\, \sum_{\{k: \, [k-1,k)\cap I_\ell \ne \emptyset\}}\chi_k T_{\la,r}h&\Bigr\|_{L^{\frac{2(n+1)}{n-1}}(M)}
	\\
	&\le r^{-\frac12}\Bigl( \, \sum_{\{k: \, [k-1,k)\cap I_\ell \ne \emptyset\}} k^{\frac{n-1}{n+1}} \|\chi_k
	T_{\la,r}h\|_{L^2(M)}^2\, \Bigr)^{\frac12} \notag 
	\\
	&\le C_Nr^{-\frac12}(1+|\la +r^{-1}(2\ell +1)|)^{\frac{n-1}{2(n+1)}}(1+|\ell|)^{-N}\|h\|_{L^2(M)},
	\notag 
\end{align*}
for each $N\in \N$, using \eqref{2.9} and the fact that if $\la_j\in [k-1,k)$ and $[k-1,k)\cap I_\ell \ne \emptyset$
then $1+r|\la_j-\la|\lesssim 1+|\ell|$ in the last step.

From this we get
\begin{align*}
	\|T_{\la,r}f\|_{L^{\frac{2(n+1)}{n-1}}(M)} &\le \sum_{\ell \in \Z}
	\bigl\|T_{\la,r}(\sum_{\la_j\in I_\ell}E_jf)\bigr\|_{L^{\frac{2(n+1)}{n-1}}(M)}
	\\
	&=\sum_{\ell\in Z}	\bigl\|\sum_{k:[k-1,k)\cap I_\ell \ne \emptyset} \chi_k T_{\la,r}(\sum_{\la_j\in I_\ell}E_jf)\bigr\|_{L^{\frac{2(n+1)}{n-1}}(M)}
	\\
	&\le C_Nr^{-\frac12}\sum_{\ell \in \Z}(1+|\ell|)^{-N}\, \bigl(1+|\la +r^{-1}(2\ell +1)|\bigr)^{\frac{n-1}{2(n+1)}}\|f\|_{L^2(M)}
	\\
	&\le Cr^{-\frac12}\la^{\frac{n-1}{2(n+1)}}\|f\|_{L^2(M)},
\end{align*}
if $N\ge 2$, since, by our assumption that $\la^{-1}r^{-1}\le 1$, we have
$$|\la +r^{-1}(2\ell +1)|=\la \bigl|1+\la^{-1}r^{-1}(2\ell +1)\bigr| \lesssim \la (1+|\ell|).$$
This concludes the proof of \eqref{2.5}.  
\end{proof}

\newsection{Localized eigenfunction estimates for other exponents}

Note that if we use \eqref{1.6}, we could repeat the proof of \eqref{2.5}
to get that for exponents $p>2$
\begin{equation}\label{3.1}
\|T_{\la,r}f\|_{L^p(M)}\le Cr^{-\frac12}\la^{\sigma(p)}\|f\|_{L^2(M)},
\quad \la \ge 1, \, \, \, \la^{-1}\le r\le \text{Inj }M,
\end{equation}
which by our earleir arguments, yields the following generalization of
\eqref{2.6}
\begin{equation*}
\|e_\la\|_{L^p(B_r(x))} \le Cr^{-\frac12}\la^{\sigma(p)}
\|e_\la\|_{L^2(B_{2r}(x))},
\quad \la \ge 1, \, \, \, \la^{-1}\le r\le \text{Inj }M.
\end{equation*}
We then could use these bounds to obtain
\begin{multline}\label{3.2}
\|e_{\la}\|_{L^p(M)}
\\
\le \la^{\sigma(p)}
\Bigl[\sup_{x\in M}r^{-\frac{p}{2(p-2)}}\|e_\la\|_{L^2(B_r)}\Bigr]^{\frac{p-2}{p}},\, \, \,
\la \ge 1, \, \, \la^{-1}\le r\le \text{Inj }M, \, \, 2<p<\infty,
\end{multline}
as well as
\begin{equation}\label{3.3}
\|e_\la\|_{L^\infty(M)}
\le C\la^{\frac{n-1}2}\sup_{x\in M}r^{-\frac12}\|e_\la \|_{L^2(B_r(x))},
\, \, \la \ge 1, \, \, \la^{-1}\le r\le \text{Inj }M.
\end{equation}
By the remarks we shall make about the relationship between these
estimates and the $L^p$-norms of the highest weight spherical harmonics,
\eqref{3.2} cannot be improved when $p_c=\tfrac{2(n+1)}{n-1}\le p<\infty$ at least for $r=\la^{-1}$, while the observations we shall make about
sup-norms of zonal functions and the right side of \eqref{3.3}
show that this estimate cannot be improved for the full range
of radii, $\la^{-1}\le r\le \text{Inj }M$.  In all liklihood
\eqref{3.2} is also sharp for this full range of $r$  if $p\ge p_c$ since
its counterpart \eqref{3.1} is best possible for this range
of $r$ and all exponents $p\ge p_c$.  The estimate \eqref{3.1} is not
sharp for $2\le p<p_c$, though.

These estimates are only of potential use for the range of exponents
$p_c\le p\le \infty$. This is because, for the
range of exponents $p_c \le p\le \infty$ the eigenfunction estimates 
\eqref{1.5} are saturated by eigenfunctions concentrating at points,
such as zonal functions on the sphere; however, for the complementary range
of exponents $2<p\le p_c$ the bounds in \eqref{1.5} are saturated
by eigenfunctions concentrating along periodic geodesics, such as
the highest weight spherical harmonics on the sphere.  We shall have
more to say about these two cases and the estimates \eqref{3.2}--\eqref{3.3} in the next section.

For the range of exponents $p_c<p\le \infty$ the author and
Zelditch in \cite{SZ1} showed that one has
\begin{equation}\label{3.4}
\|e_\la\|_{L^p(M)}=o(\la^{\sigma(p)})
\end{equation}
for generic $(M,g)$, and in recent papers \cite{SZ2}--\cite{SZ3} gave
necessary and sufficient conditions in the real analytic case for
a stronger version of \eqref{3.4} involving quasimodes.

For the complementary range of exponents $2<p<p_c$, because of reasons that
we just alluded to, one would not expect bounds like \eqref{3.3} to be useful for proving \eqref{3.4}.  Instead, in a series of papers
of the author \cite{SK} and Blair and the author \cite{BS1}--\cite{BS2},
motivated by earlier related work of Bourgain~\cite{Bourgain},
the strategy was to prove a variation of \eqref{3.2} which controls
the $L^p$-norms of eigenfunctions in terms of their $L^2$-mass on
small tubes about geodesics.  Specifically if $\varPi$ denotes the
space of unit length geodesics in $M$ and if 
${\mathcal T}_{\la^{-\frac12}}(\gamma)$ denotes a $\la^{-\frac12}$-tube about a given $\gamma\in \varPi$, it was shown that
\begin{equation}\label{3.5}
\|e_\la\|_{L^p(M)}\le C\la^{\sigma(p,n)}
\, \Bigl[ \, \sup_{\gamma\in \varPi}
\int_{{\mathcal T}_{\la^{-\frac12}}(\gamma)}|e_\la|^2 \, dV_g
\, \Bigr]^{\theta(p,n)}, \, \, \,
\la \ge 1, \, \, 2<p<p_c,
\end{equation}
for some 
$$\theta(p,n)>0.$$

When $n=2$ the author and Zelditch \cite{SZStein} were able to show that 
if $(M,g)$ has nonpositive curvature one has
\begin{equation}\label{3.6}
\sup_{\gamma\in \varPi}\int_{{\mathcal T}_{\la^{-\frac12}}(\gamma)}|e_\la|^2 \, dV_g=o(1),
\end{equation}
and, therefore, by \eqref{3.5}, one gets the improved eigenfunction bounds
\eqref{3.4} when $n=2$ if $2<p<p_c$.  For higher dimensions, Blair
and the author~\cite{BS2} obtained \eqref{3.6} and hence 
\eqref{3.4} under this curvature assumption.

For the elusive endpoint case $p=p_c$ Hezari and Rivi\`ere~\cite{HR} were
able to obtain a stronger version of 
\eqref{3.4} involving logarithmic improvements for a density one
subsequence of eigenfunctions on a given manifold of negative curvature.
They did this by proving results like Theorem~\ref{loctheorem}
when $r$ is a power of $1/\log \la$
and then obtaining, for a density one sequence of eigenfunctions, very
natural $L^2(B_r)$-norms for such $r$.  We shall say more about
the latter results, which were also obtained independently
by Han~\cite{H}, in the next section.

\newsection{Remarks on the size of $L^2(B_r)$}

Let us conclude with a few remarks about the size of
$L^2(B_r(x))$-norms of eigenfunctions.  

The first is that for any
$(M,g)$ we have the trivial estimates
\begin{equation}\label{4.1}
\|e_\la\|_{L^2(B_r(x))}\le Cr^{\frac12}, 
\quad x\in M, \, \, \la \ge 1, \, \, \la^{-1}\le r\le \text{Inj }M.
\end{equation}
An interesting question would be to determine when one
can improve on this easy estimate for $r=r(\la)\to 0$
as $\la \to \infty$ (either through all eigenvalues
or subsequences)
using some dynamical or geometric assumption, such as $(M,g)$ having everywhere
nonpositive curvature.  We shall go over a few model cases after
presenting the proof of \eqref{4.1}.

To prove this inequality,
we may of course assume that $\text{Inj }M\ge 4$.
Then,
if $\rho\in {\mathcal S}(\R)$ is as in \eqref{2.1},
then it suffices to show that
$$
\|\rho(\la-P)f\|_{L^2(B_r(x))}\le Cr^{\frac12}\|f\|_{L^2(M)},
\quad x\in M, \, \, \la \ge 1, \, \, \la^{-1}\le r\le \text{Inj }M,
$$
since $e_\la =\rho(\la-P)e_\la$.
By a routine $TT^*$ argument, if $\eta =|\rho|^2$, this is equivalent
to 
\begin{multline}\label{4.2}
\|\eta(\la-P)h\|_{L^2(B_r(x))}\le Cr\|h\|_{L^2(B_r(x))},
\\
\text{if } \, \text{supp }h\subset B_r(x), \, \, \text{and } \, 
  \la \ge 1, \, \, \la^{-1}\le r\le \text{Inj }M.
\end{multline}
If we argue as in the proof of \cite[Lemma 5.1.2]{SBook} we find
that, since $\supp \Hat \eta \subset [-2,2]$ and we are assuming that
$\text{Inj }M\ge 4$, the kernel of $\eta(\la-P)$ can be written as
$$
\la^{\frac{n-1}2} \sum_\pm a_\pm\bigl(\la, d_g(x,y)\bigr)
\, \bigl(d_g(x,y)\bigr)^{-\frac{n-1}2}e^{\pm i\la d_g(x,y)},
\quad \text{if } \, d_g(x,y)\ge \la^{-1},
$$
where
$$
\bigl|\partial_s^ja_\pm(\la,s)\bigr|\le C_j s^{-j}, 
\quad s\ge \la^{-1}, \, \, \text{and } \, j=0,1,2,\dots,
$$
and this kernel is $O(\la^{n-1})$ when $d_g(x,y)\le \la^{-1}$.
Using this, it is routine to obtain \eqref{4.2} using
H\"ormander's $L^2$-oscillatory integral estimates
(e.g. \cite[Theorem 2.1.1]{SBook}).  The argument one uses
is very similar to the proof of \cite[Theorem 5.2.1]{SBook}.

As one would expect, the estimates \eqref{4.1} are saturated
on the standard spheres.  In that case the $L^2$-normalized zonal
eigenfunctions $Z_\la$ centered at a given
point $x_0\in S^n$ with $\la =\sqrt{k(k+n-1)}$ are given
by the formula
$$
Z_\la(x) = (d_k|S^n|)^{-\frac12}\sum_{j=1}^{d_k}e_{k,j}(x)
\overline{e_{k,j}(x_0)}
$$
where $d_k\approx \la^{n-1}$ denotes
the dimension of the space of spherical harmonics of degree
$k$ and $\{e_{k,j}\}_{j=1}^{d_k}$ is an orthonormal basis
of this space (see e.g., \cite[\S 3.4]{SoHang} for more details).
One can use the classical Darboux formula (see e.g. \cite[4.7]{Sthesis})
to see that 
\begin{multline*}
Z_\la(x)\approx \bigl(d_g(x,x_0)\bigr)^{-\frac{n-1}2}
\Bigl[\cos\bigl(N_kd_g(x,x_0)+\gamma\bigr) \, +\, O(1)(\la d_g(x,x_0))^{-1}
\Bigr],
\\
\text{if }\, \, N_k=(2k+n-1)/2, \, \,
 \gamma=-(n-1)\pi/4, \, \, \, \text{and } \, \, d_g(x,x_0)\ge \la^{-1},
\end{multline*}
as well as
$$
Z_\la(x,x_0)=O(\la^{\frac{n-1}2}) \quad
\text{if } \, \, d_g(x,x_0)\le \la^{-1}.
$$
Using these facts, we find that
$$
\|Z_\la\|_{L^2(B_r(x_0))} \approx r^{\frac12},
\quad \text{if } \, \, \la^{-1}\le r\le \pi,
$$
showing, as claimed, that \eqref{4.1} is saturated on $S^n$.

For the zonal functions, we have just shown that we get
{\em no improvement} for the size of
$r^{-\frac12}\|e_\la\|_{L^2(B_r)}$ by taking $r$ to be very
small.  For the other extreme spherical harmonics, the
highest weight spherical harmonics, $Q_\la$, as it turns out,
one does achieve an improvement by taking $r$ to be small.
Recall that the $Q_\la$ are the restriction of the harmonic
polynomials $c_k(x_1+ix_2)^k$ to the unit sphere $|(x_1,x_2, \dots, x_{n+1})|=1,$
and $c_k\approx k^{\frac{n-1}4}$ for $L^2$-normalization.  Therefore,
$$
|Q_\la| \approx k^{\frac{n-1}4} e^{-\frac{k}2 (x_3^2+\dots+x_{n+1}^2)}
$$
and $Q_\la$ is an eigenfunction with frequency $\sqrt{k(k+n-1)}=\la$,
as above.  Thus, $Q_\la$ behaves like a Gaussian beam concentrated
along a $\la^{-\frac12}$ neighborhood of the set on the unit
sphere centered at the origin in ${\mathbb R}^{n+1}$ where $x_3= \dots =x_{n+1}=0$,
which, of course is a periodic geodesic.  Using these size
estimates we find that if $y$ is a point on this geodesic then
$$
\|Q_\la\|_{L^2(B_r(y))} \approx r^{\frac12}, \quad
\la^{-\frac12}\le r\le \pi
$$
meaning that there is no improvement over \eqref{4.1} for this range
of $r$, while on the other hand,
we have the uniform bounds as $x$ ranges over $S^n$:
\begin{equation}\label{4.3}
r^{-\frac12}\|Q_\la\|_{L^2(B_r(x))}\le Cr^{\frac{n-1}4}, \quad
\text{if } \, r=\la^{-1}.
\end{equation}
Curiously, if we use this fact along with \eqref{3.2}--\eqref{3.3}
we get
$$
\|Q_\la\|_{L^p(S^n)}\lesssim \la^{\frac{n-1}2(\frac12-\frac1p)},
\quad \tfrac{2(n-1)}{n-1}\le p\le \infty,
$$
which is sharp since, in view of the above $\|Q_\la\|_{L^p(S^n)}$
is comparable to $\la^{\frac{n-1}2(\frac12-\frac1p)}$ for
all $p\ge2$.  One also gets smaller improvements on these $L^2(B_r)$-norms
for $Q_\la$ if $\la^{-1}< r\ll \la^{-\frac12}$.

In the case of manifolds of nonpositive curvature, we can also beat the
trivial estimate \eqref{4.1} if we use a result of B\'erard~\cite{PB}
if $(M,g)$ is of nonpositive curvature.  Recall that, in this case,
he showed that the error term in the local Weyl law is $O(\la^{n-1}/\log \la)$,
which implies that
\begin{equation}\label{4.4}
|e_\la(x)|\le C\la^{\frac{n-1}2}/(\log \la)^{\frac12}.
\end{equation}
Using this, we deduce that we have
\begin{equation}\label{4.5}
\|e_\la\|_{L^2(B_r(x))}\le C\frac{r^{\frac12}}{(\log \la)^{\frac12}},
\quad \text{if } \, \, r=\la^{-1},
\end{equation}
under this curvature assumption.  The argument is circular, but we can
then use \eqref{3.3} to obtain \eqref{4.4}.
Unfortunately we cannot use \eqref{4.5} along with \eqref{3.2} to get 
any improvements over \eqref{1.4} if $p<\infty$ for manifolds of nonpositive curvature,
although, it was already known 
by a recent result of Hassell and Tacy~\cite{HT}
that for $\tfrac{2(n+1)}{n-1}<p<\infty$, like in \eqref{4.4}, one gets
a $(\log \la)^{-\frac12}$ improvement over \eqref{1.4} in this case.

If one works on the much larger scale $r=1/(\log \la)^{\kappa}$,
where the power $\kappa=\kappa_n$ depends on $n$, Han~\cite{H}
and Hezari and Rivi\`ere~\cite{HR} showed that there is a density one
sequence of eigenvalues, $\{\la_{j_k}\}$ so that
\begin{equation}\label{4.6}
\int_{B_r(x)} |e_{\la_{j_k}}|^2 \, dV_g=
\frac{|B_r(x)|}{|M|}+o(r^n), \quad r=1/(\log \la)^\kappa
.
\end{equation}
Since
the radii shrink as $\la$ increases, this of course does not follow from
the classical quantum ergodic identity \eqref{1.10}.  The latter
holds whenever the geodesic flow is ergodic, while special dynamical properties
of the geodesic flow on negatively curved
manifolds was used in the aforementioned results
to obtain \eqref{4.6}.  Using \eqref{1.2} one can obtain certain
$\log$-power improvements of \eqref{1.3} for density one subsequences
of eigenfunctions on negatively curved manifolds as was done in
Hezari and Rivi\'ere~\cite{HT}.

Han asked in \cite{H} whether one could break the logarithmic barrier
and prove the variant of \eqref{4.6} where
$r=\la^{-\kappa}$ for some $0<\kappa<1$.  An affirmative answer to this
seemingly difficult question would similarly lead to $\la^{-\kappa
\frac{n-1}{2(n+1)}}$ improvements of \eqref{1.3}.  We remark, though,
that one does not need the full strength of \eqref{4.6} to get improvements
over \eqref{1.3}.  Indeed, if one could show that
\begin{equation}\label{4.7}
\sup_{x\in M}\|e_\la\|_{L^2(B_{r(\la)}(x))}\le C\bigl(r(\la)\bigr)^{\frac{n}2},
\quad r(\la)\to 0,
\end{equation}
as $\la$ ranges over a subsequence $\la_{j_k}$ of eigenvalues, then
Theorem~\ref{loctheorem} would yield
\begin{equation}\label{4.8}
\|e_\la\|_{L^{\frac{2(n+1)}{n-1}}(M)}\le C\bigl(r(\la) \, \la\bigr)^{\frac{n-1}{2(n+1)}}
\end{equation}
for this subsequence.
To obtain the missing $1/(\log \la)^{1/2}$ 
endpoint result with $p=p_c$ of Hassell and Tacy~\cite{HT}
that we mentioned before would require $r(\la)=1/(\log \la)^{\frac{n+1}{n-1}}$, which involves a power, $\tfrac{n+1}{n-1}>1$, which is larger
than the ones occurring in \cite{H} and \cite{HR}.

\subsection*{Acknowledgements}

The author is grateful to Hamid Hezari and Gabriel Rivi\`ere for going through an early draft of this paper, as well as for a very informative talk
by Hezari at Johns Hopkins University explaining his joint work with Rivi\`ere.

\end{document}